\documentclass[notitlepage,11pt]{article}
\usepackage{amssymb,mathrsfs,amsmath}
\catcode`\@=11
\@addtoreset{equation}{section}

\catcode`\@=12
\usepackage{amsmath}
\usepackage{latexsym}
\usepackage{amscd}
\usepackage{xspace}
\usepackage{verbatim}

\newtheorem{Theorem}{Theorem}[section]
\newtheorem{Lemma}[Theorem]{Lemma}

\newtheorem{Remark}[Theorem]{Remark}

\newcommand{\Rn}{{\mathbb R}^{N}}

\newcommand{\ba} {\beta}

\newcommand{\De} {\Delta}
\newcommand{\la} {\lambda}

\newcommand{\var} {\varepsilon}

\newcommand{\R}{\mathbb{R}}
\def\proof{\noindent{\textbf{Proof. }}}
\def\QED{\hfill {$\square$}\goodbreak \medskip}

\newcommand{\un}{u_{n}}
\newcommand{\deb}{\rightharpoonup}
\newcommand{\starstar}{p^{*\!*} }

\begin{document}

\title{Entire solutions for a class of elliptic equations involving p-biharmonic operator and Rellich potentials}

\author{Mousomi Bhakta
\footnote
{Department of Science and Technology, University of New England, Armidale, NSW-2350, Australia. 
Email: {mousomi.bhakta@gmail.com}} 
}

\date{}

\maketitle

\begin{abstract}
\noindent
\footnotesize We study existence, multiplicity and qualitative properties of entire
solutions for a noncompact problem related to p-biharmonic type equations with weights. More precisely, we deal with the following family of equations
$$
{\Delta_{p}^2u=\lambda|x|^{-2p}|u|^{p-2}u+|x|^{-\beta}|u|^{q-2}u\quad\textrm{in $\R^N$,}}
$$
where $N> 2p$, 
$p>1$, $q>p$, $\beta\!=\!N\!-\frac{q}{p}(N-2p)$ and $\lambda\!\in\!\R$ is smaller than the Rellich constant.

%
\bigskip

\medskip

\end{abstract}

\section{Introduction}

In this article  we  study  weak solutions to the problem
\begin{equation}
\label{eq:problem}
\begin{cases}
\Delta^2_p u=\lambda|x|^{-2p}|u|^{p-2}u+|x|^{-\beta}|u|^{q-2}u&\textrm{in $\R^N$}\\
u\in  D^{2,p}(\R^N)~,\quad u\neq 0~,
\end{cases}
\end{equation}
where
$\Delta_p^2 u=\Delta(|\Delta u|^{p-2}\Delta u)$ and
\begin{equation}
\label{eq:ass_q_lambda}
p>1, N>2p, \quad q>p, \quad \beta=N-\frac{q}{p}(N-2p), \quad \lambda<\gamma_{N,p}^p \  \text{and}\  \gamma_{N,p}=\frac{N(p-1)(N-2p)}{p^2}.
\end{equation}
Here $D^{2,p}(\Rn)$ is the closure of $C^\infty_0(\Rn)$ with respect to the norm $(\int_{\Rn}|\De u|^p dx)^\frac{1}{p}$.
It is well known that $\gamma_{N,p}^p$ is the
best constant in the  Rellich inequality
\begin{equation}\label{Rellich}
\gamma_{N,p}^p\int_{\Rn}|x|^{-2p}|u|^p dx\displaystyle\leq \int_{\Rn}|\De u|^p dx  \quad\textrm{for any}\   \  u\in  D^{2,p}(\R^N).
\end{equation}
 In literature \eqref{Rellich} with $p=2$ is considered as the classical Rellich inequality and was proved by Rellich in 1953 (See \cite{Rel54}, \cite{Rel69}). Later Davies and Hinz \cite{DH} generalized the classical Rellich inequality and showed that \eqref{Rellich} holds for any $p\in (1,\frac{N}{2})$.

If $q$ coincides with the critical Sobolev exponent
$$
\starstar:=\frac{Np}{N-2p}~\!,
$$
then (\ref{eq:problem}) becomes
\begin{equation}
\label{eq:problem_critical}
\begin{cases}
\Delta^2_p u=\lambda|x|^{-2p}|u|^{p-2}u+|u|^{\starstar-2}u&\textrm{in $\R^N$}\\
u\in  D^{2,p}(\R^N)~,\quad u\neq 0.
\end{cases}
\end{equation}
When $\la=0$, it's well known from the celebrated paper of P.L.Lions \cite{Li84} that  \eqref{eq:problem_critical} has a  positive solution $U$ which is the extremal for the Sobolev inequality 
\begin{equation}\label{Sobolev}
S^{**}(\int_{\Rn}|u|^{\starstar} dx)^\frac{p}{\starstar}\leq \int_{\Rn}|\De u|^p dx
\end{equation}
where $S^{**}$ is the Sobolev Constant.
To prove that $S^{**}$ is achieved, P.L. Lions had shown that every bounded minimizing sequence is relatively compact up to dilations and translations. Moreover, by using Schwarz symmetrization he showed that, up to a change of sign, any extremal for $S^{**}$ is radially symmetric, nonnegative and decreasing. Later using this information, Hulshof and Van der Vorst \cite{HV} proved the uniqueness of extremals for $S^{**}$, modulo dilations, translations in $\Rn$ and change of sign.

If $p\leq q\leq\starstar$ and $\beta$ are as in (\ref{eq:ass_q_lambda}), then by interpolating  (\ref{Rellich}) and (\ref{Sobolev}) via H\"older inequality, it can be easily shown that there exists a constant $C=C(N,p, q)>0$ such that
\begin{equation}\label{CKN}
 C\left(\int_{\Rn}
|x|^{-\beta}|u|^{q}dx\right)^{p/q}\leq \int_{\Rn}|\De u|^{p}dx \quad \textrm{for any $u\in D^{2,p}(\Rn)$}.
\end{equation}

Notice that  (\ref{CKN}) with $p=2$ is the fourth-order version of the celebrated Caffarelli-Kohn-Nirenberg inequalities \cite{CKN}. We cite also \cite{CM2} for a large class of dilation-invariant inequalities on cones. 

In recent years problems related with the inequality (\ref{CKN}) (in the case $p=2$) and the equation with biharmonic operator have been  investigated in several works, we quote \cite{AlDo} \cite{BM}, \cite{C},  \cite{CM2}, \cite{CatWan},  \cite{Cof}, \cite{GG},  \cite{Miti},\cite{Mu1}, \cite{Mu2} and the references there-in. Recently the generalized version of the inequality \eqref{CKN} and the extremal of that inequality has been studied by Roberta Musina. (see  \cite{Mu1} and \cite{Mu2}). \\

Note that that the choice of $\ba$ in \eqref{eq:ass_q_lambda} makes \eqref{eq:problem}  invariant with respect to the weighted dilation 
\begin{equation}\label{transform}
u(x) \mapsto t^\frac{N-2p}{p}u(tx), \quad (t>0) .
\end{equation}
As a consequence, the corresponding variational problems exhibit a lack of compactness. 

It is clear that the infimum
$$
S_q(\la):  =  
\inf_{\scriptstyle u\in  D^{2,p}(\Rn)
\atop\scriptstyle u\ne 0}
\frac{\displaystyle
\int_{\Rn}|\De u|^{p}dx-\la\int_{\Rn}|x|^{-2p}|u|^{p}dx}
{\displaystyle\left(\int_{\Rn}|x|^{-\beta}|u|^{q}dx\right)^{p/q}}
$$
is positive, provided that $\lambda<\gamma_{N,p}^p$.
In addition, extremals for $S_q(\la)$ give rise to solutions to (\ref{eq:problem}) upto a multiplicative constant.

Define,
$$
S_q^{\rm rad}(\la):  =  \inf_{\scriptstyle u\in  D^{2,p}(\Rn)
\atop\scriptstyle u=u(|x|)~,~u\ne 0}\frac{\displaystyle
\int_{\Rn}|\De u|^{p}dx-\lambda\int_{\Rn}|x|^{-2p}|u|^{p}dx}
{\displaystyle\left(\int_{\Rn}|x|^{-\beta}|u|^{q}dx\right)^{p/q}}
$$ which is positive as $S_q^{\rm rad}(\la)\geq S_q(\la)>0$ when $\la<\gamma_{N,p}^p$. 

\vspace{2mm}

We prove the following results.  

\begin{Theorem}
\label{T:existence_radial}
Let $N>2p$, $q>p$, $\lambda<\gamma_{N,p}^p$, and 
$\beta=N-\frac{q}{p}(N-2p)$. Then the problem (\ref{eq:problem}) has at least one radially symmetric solution $u$ which achieves $S_q^{\rm rad}(\la)$.
\end{Theorem}
In fact when $p=2$ and $-(N-2)^2\leq\la<(\gamma_{N,2})^2$, problem \eqref{eq:problem} has a unique radial solution. (See \cite{BM}).

For $\la=0$, Theorem \ref{T:existence_radial} was proved in \cite[Theorem 1.3]{Mu2}. Following the same procedure as in \cite{Mu2}, Theorem \ref{T:existence_radial} can be proved in the case $\la\not=0$ as well.

\medskip

In case $q\le \starstar$ one can use again variational methods to find
solutions to (\ref{eq:problem}) that are not necessarily radially symmetric.

In Section \ref{S:non_radial} we
prove the next existence result.

\begin{Theorem}\label{exnonex}
Let $N>2p$, $q\in(p,\starstar]$, $\lambda<\gamma_{N,p}^p$, and 
$\beta=N-\frac{q}{p}(N-2p)$. Then
\begin{enumerate}
\item [(i)]
The infimum $S_q(\la)$ is achieved for any $q\in(p,\starstar)$.
\item[(ii)]
The infimum $S_{\starstar}(\la)$ is achieved if and only if
$\la\ge 0$.
\end{enumerate}
\end{Theorem}

In Section \ref{S:symmetry} we wonder whether the solutions in Theorem
\ref{exnonex} are radially symmetric or breaking symmetry occurs. First, by using rearrangement techniques
we prove that extremal for $S_q(\la)$ is always radially symmetric provided that $\lambda\ge 0$.
In contrast, we show that if {\bf $\lambda<<0$} and {\bf $p\geq 2$} then 
{\bf $S_q(\lambda)<S_q^{\rm rad}(\la)$} i.e break of symmetry occurs. Therefore, if in addition $q\in (p, \starstar)$ then problem \eqref{eq:problem} has at least two distinct nontrivial solutions.

\section{Existence and non existence of ground state solutions}
\label{S:non_radial}
In this section we will prove Theorem \ref{exnonex}. A key tool in our argument is the following $\varepsilon-\text{compactness}$ lemma. This result is an adaptation of a tool already used in previous works, like \cite{BM} or \cite{CM2}. Therefore we omit the proof.

\begin{Lemma}\label{lemma1}
Let $\un\in D^{2,p}(\Rn)$ such that $\un\deb 0$ in $\mathcal D^{2,p}(\Rn)$ and 
\begin{equation}\label{assump1}
\De^2_p\un-\la |x|^{-2p}|\un|^{p-2}u_n=|x|^{-\ba}|\un|^{q-2}\un+f_n
\end{equation}
\begin{equation}\label{assump2}
\int_{B_{R}}|x|^{-\ba}|\un|^q dx\leq \var_{0}\quad\textit{for some $\var_0, R>0$,}
\end{equation}
where $f_n\to 0$ in  
the dual space of $ D^{2,p}(\Rn)$. If $\var_{0}<S_{q}(\la)^{\frac{q}{q-p}}$,
then 
$$
|x|^{-\ba}|\un|^q\to 0\quad\textit{in $L^1_{\rm loc}(B_R)$.}
$$
\end{Lemma}

\subsection{Proof of Theorem \ref{exnonex} (i)}

{\bf Step 1:} Using Ekeland's variational principle we can choose a minimizing sequence $\{u_n\}$ for
$S_q(\lambda)$ such that
\begin{eqnarray}\label{assump4}
\int_{\Rn}|\De\un|^p-\la\int_{\Rn}{|x|^{-2p}}|\un|^p &=& \int_{\Rn}|x|^{-\ba}|\un|^q dx\nonumber
\\
&=& S_{q}(\la)^\frac{q}{q-p}+o(1)~\!,
\end{eqnarray}
\begin{equation}\label{assump3}
\Delta^2_p\un-\la |x|^{-2p}|\un|^{p-2}\un=|x|^{-\ba}|\un|^{q-2}\un+f_n~\!,
\end{equation}
where $f_n\to 0$ in the dual space of $ D^{2,p}(\Rn)$. Up to a rescaling, we  assume that
\begin{equation}\label{a}
\int_{B_{2}}|x|^{-\ba}|\un|^q dx=\frac{1}{2}S_{q}(\la)^\frac{q}{q-p}~\!.
\end{equation} 
Therefore it can be easily checked that $\un$ is a bounded sequence in
$D^{2,p}(\Rn)$ by Rellich inequality as we have $\la<\gamma_{N,p}^p$. Hence we can assume that 
there exists $u\in  D^{2,p}(\Rn)$ such that $\un\deb u$ 
weakly in
$ D^{2,p}(\Rn)$. \\

{\bf Claim 1:} $u\neq 0$.\\

We argue by contradiction. Suppose $\un\deb 0$. Then using  Lemma \eqref{lemma1} we obtain
$$
o(1)=\int_{B_{1}}|x|^{-\ba}|\un|^q dx=\int_{B_{2}}|x|^{-\ba}|\un|^q dx-
\int_{1<|x|<2}|x|^{-\ba}|\un|^q dx~\!.
$$
Thus from  (\ref{a}) we infer
\begin{equation}
\label{eq:new}
\int_{1<|x|<2}|x|^{-\ba}|\un|^q dx=\frac{1}{2}S_{q}(\la)^\frac{q}{q-p}+o(1)
\end{equation}
which leads to a contradiction by Rellich's compactness theorem, as $q\in(p,\starstar)$. Thus Claim 1 follows.

\vspace{2mm}
 
{\bf Claim 2:} $u$ is a weak solution of \eqref{eq:problem}.

For $p=2$, it is straight forward. For $p\not=2$, we choose $w_n\in D^{2,p}(\Rn)$ such that 
$$\Delta(|\Delta w_n|^{p-2}\Delta w_n)=f_n; \quad\text{and}\quad w_n\to 0 \quad\text{in}\quad  D^{2,p}(\Rn),$$
where $f_n$ is as in \eqref{assump3}. Define, 
$h_n :=-|\Delta u_n|^{p-2}\Delta u_n+|\Delta w_n|^{p-2}\Delta w_n$, that is bounded in $L^{\frac{p}{p-1}}(\Rn)$. Also from \eqref{assump3}, it follows that,
\begin{equation}\label{eq:h_n}
-\Delta h_n=\la|x|^{-2p}|u_n|^{p-2}u_n+|x|^{-\ba}|u_n|^{q-2}u_n.
\end{equation}
As $p<q<p^{**}$, it is easy to see that 1st term of RHS of \eqref{eq:h_n} is bounded in $L^\frac{q}{p-1}_{loc}(\Rn\setminus\{0\})$ and 2nd term of RHS of \eqref{eq:h_n} is bounded in $L^\frac{q}{q-1}_{loc}(\Rn\setminus\{0\})$. Since $p<q$ implies $\frac{q}{q-1}<\frac{q}{p-1}$, we can easily conclude that $-\Delta h_n$ is bounded in $L^\frac{q}{q-1}_{loc}(\Rn\setminus\{0\})$, which in turn implies $h_n$ is bounded in $W^{2,\frac{q}{q-1}}_{loc}(\Rn\setminus\{0\})$. Therefore almost everywhere convergence of  $h_n$ follows and so of $u_n$, as $w_n\to 0$ a.e. Hence using Vitaly's convergence theorem via Holder inequality, we obtain
$$\int_{\Rn}|\Delta u_n|^{p-2}\Delta u_n\Delta\phi\to\int_{\Rn}|\Delta u|^{p-2}\Delta u\Delta\phi \quad \forall\quad \phi\in C^{\infty}_0(\Rn).$$ 
Therefore it is easy to see that the claim follows.

\vspace{2mm}

{\bf Step 2:} Therefore $S_q(\la)\displaystyle\left(\int_{\Rn}\frac{|u|^q}{|x|^{\ba}}dx\right)^\frac{p}{q}\leq\int_{\Rn}\left(|\Delta u|^{p}-\la\frac{|u|^p}{|x|^{2p}}\right)dx=\int_{\Rn}\frac{|u|^q}{|x|^{\ba}}dx $. \\
Since $u\not= 0$, we obtain $\displaystyle\int_{\Rn}\frac{|u|^q}{|x|^{\ba}}dx\geq S_q(\la)^\frac{q}{q-p}$. On the other hand by the lower semicontinuity of the norm in $L^q(\Rn, |x|^{-\ba}dx)$ and \eqref{assump4}, we have 
$\displaystyle\int_{\Rn}\frac{|u|^q}{|x|^{\ba}}dx\leq S_q(\la)^\frac{q}{q-p}$. \\
Therefore $\displaystyle\int_{\Rn}\frac{|u|^q}{|x|^{\ba}}dx= S_q(\la)^\frac{q}{q-p}$ which in turn implies that $S_q(\la)$ is achieved by $u$.

\QED

\subsection{Proof of Theorem \ref{exnonex} in the limiting case $q=\starstar$}
We start by pointing out a sufficient condition for existence.

\begin{Lemma}\label{prop2}
If $S_{\starstar}(\la)<S^{*\!*}$ then $S_{\starstar}(\la)$ is achieved.
\end{Lemma}

\proof
As in the proof of Theorem \ref{exnonex} $(i)$, we choose a minimizing sequence $\{u_n\}$ for
$S_{\starstar}(\la)$ satisfying
\begin{eqnarray*}
\int_{\Rn}|\De\un|^p~\!dx-\la\int_{\R^N}{|x|^{-2p}}|\un|^p~\!dx&=& \int_{\Rn}|\un|^{\starstar} dx+o(1)\\
&=&
 S_{\starstar}(\la)^{N/2p}+o(1)
 \end{eqnarray*}
 \begin{equation}\label{assump6}
\Delta^2_p\un-\la |x|^{-2p}|\un|^{p-2}\un=|\un|^{\starstar-2}\un+f_n
\end{equation}
$$
\int_{B_{2}}|\un|^{\starstar} dx = \frac{1}{2} S_{\starstar}(\la)^{N/2p}~\!,
$$
where $f_n\to 0$ in the dual of $ D^{2,p}(\Rn)$. In addition, we can assume that
$\un$ weakly converges to $u$  in $D^{2,p}(\Rn)$. 

{\bf Claim 1:} $u\ne 0$.\\

 We will prove the claim by contradiction, thus we assume 
$\un\deb 0$ in $D^{2,p}(\Rn)$. 
Arguing as in the proof of Theorem \ref{exnonex} $(i)$ we can conclude that (\ref{eq:new}) holds with
$q=\starstar$, that is,
\begin{equation}\label{b}
\int_{1<|x|<2}|\un|^{\starstar}~\!dx=\frac{1}{2}S_{\starstar}(\la)^{N/2p}+o(1).
\end{equation}
Now we choose  a cut off function  $\phi$ in $C_{c}^{\infty}(\Rn\setminus\{0\})$ such that $\phi\equiv 1$ in $B_2\setminus B_1$. Taking  $\phi^p\un$ as a test function in \eqref{assump6} and using Rellich's compactness theorem and H\"older inequality we obtain
$$
\int_{\Rn}[|\De(\phi\un)|^p-\la|x|^{-2p}|\phi\un|^p]dx \leq S_{\starstar}(\la)\displaystyle\left(\int_{\Rn}|\phi\un|^{\starstar}dx\right)^\frac{p}{\starstar}+o(1)
$$
 Since $\phi$ has a compact support in
$\R^N\setminus\{0\}$,  using Rellich's compactness theorem and the Sobolev inequality and arguing as in the proof of Theorem \ref{exnonex} $(i)$ we obtain 
$$
\int_{\Rn}[|\De(\phi\un)|^p-\la|x|^{-2p}|\phi\un|^p]=
\int_{\Rn}|\De(\phi\un)|^p+o(1) \ge
S^{*\!*}\displaystyle\left(\int_{\Rn}|\phi\un|^{\starstar}\right)^{p/\starstar}.
$$
Therefore we have, 
$$
S^{*\!*}\displaystyle\left(\int_{\Rn}|\phi\un|^{\starstar}dx\right)^{p/\starstar}\leq S_{\starstar}(\la)\displaystyle
\left(\int_{\Rn} |\phi\un|^{\starstar}dx\right)^{p/\starstar}~\!.
$$ 
As $S_{\starstar}(\la)<S^{**}$, the above inequality
 implies $\int_{\Rn}|\phi\un|^{\starstar} =o(1)$.
Hence $\int_{B_2\setminus B_1}|\un|^{\starstar}dx=o(1)$, since
$\phi\equiv 1$ in $B_2\setminus B_1$.  This is a contradiction to \eqref{b}. Thus $u\not=0$. 

It follows by a standard concentration-compactness technique by P. L. Lions (see the proof of \cite[Theorem 2.4]{Li84}) that $u_n$ is relatively compact  and therefore up to a subsequence $u_n\to u$ in $D^{2,p}(\Rn)$. Hence $S_{\starstar}(\la)$ is achieved.
 
\QED

\medskip
\noindent
{\bf Proof of Theorem \ref{exnonex} - $\bf (ii)$.} First we show that $S_{\starstar}(\la)<S^{**}$ holds when $\la\in (0,\gamma_{N,p}^p)$. Let $U$ be an extremal of $S^{**}$ which exists by \cite{Li84}. Therefore if $\la\in (0,\gamma_{N,p}^p)$, then we have
$$
S_{\starstar}(\la) \leq\frac{\displaystyle\int_{\Rn} \displaystyle\left[|\De U|^p-\la|x|^{-2p}|U|^p\right]~\!}{\displaystyle\left(\displaystyle\int_{\Rn}|U|^{\starstar}~\!\right)^{p/\starstar}}
< \frac{\displaystyle\int_{\Rn}|\De U|^p~\!}{\left(\displaystyle\int_{\Rn}|U|^{\starstar}~\!\right)^{p/\starstar}}= S^{*\!*}.
$$
 hence $S_{\starstar}(\la)$ is achieved.\\
 Now it remains to study the case $\la<0$. In this case, it is easy to see that $S_{\starstar}(\la)\geq S^{**}$. Now we choose an arbitrary function
$u$ in $C^\infty_c(\R^N\setminus\{0\})$ and  set
$u_y(x)=u(x+y)$. Therefore
\begin{eqnarray*}
S_{\starstar}(\la) &\leq& \lim_{|y|\to\infty}\frac{\displaystyle\int [|\De u_y|^p-\la|x|^{-2p}|u_y|^p]~\!}{\displaystyle\left(\displaystyle\int_{\Rn}|u_y|^{\starstar}~\!\right)^{p/\starstar}}
= \lim_{|y|\to\infty}\frac{\displaystyle\int [|\De u|^p-\la|x-y|^{-2p}|u|^p]~\!}{\displaystyle\left(\displaystyle\int_{\Rn}|u|^{\starstar}~\!\right)^{p/\starstar}}
\\
&=&
\frac{\displaystyle\int |\De u|^p~\!}{\displaystyle\left(\displaystyle\int_{\Rn}|u|^{\starstar}~\!\right)^{p/\starstar}},
\end{eqnarray*}
which implies $S_{\starstar}(\la)\leq S^{**}$. Hence $S_{\starstar}(\la)= S^{*\!*}$. Therefore $S_{\starstar}(\la)$ can not be achieved since $S^{*\!*}$ is achieved.

\QED

\section{Nonnegativity, symmetry and breaking symmetry}
\label{S:symmetry}

In this section we will
study the symmetry, nonnegativity and breaking symmetry of the extremal of $S_{q}(\la)$ depending on the parameter $\la$. It's known from \cite{Li84} that when $\la=0$ and $q =\starstar$, the Sobolev constant $S^{**}$ is achieved by a radially symmetric nonnegative and decreasing function.

Since truncations $u\mapsto u^{\pm}$ are not allowed in dealing with fourth order differential operators, the nonnegativity of extremals for $S_q(\la)$ does not follow by usual arguments. 

\begin{Theorem}
\label{T:symmetry}
Assume $\lambda\neq 0$ or $q<\starstar$.
If $\la\geq 0$ then $S_{q}(\la)$ is achieved by a positive function $u\in D^{2,2}(\R^N)$. Moreover, $u$ is radially symmetric about the origin and radially decreasing.
\end{Theorem}
The proof is based on rearrangement technique which was already used to prove this result in the case $p=2$ in \cite{BM}. This is an easy adaptation 
of the proof used for $p=2$.

\proof
Let $u$ be an extremal of $S_{q}(\la)$ and we 
denote by  $(-\De u)^{*}$ the Schwarz symmetrization of $-\De u$. Let
 $v\in D^{2,p}(\Rn)$ such that
$$-\De v=(-\De u)^{*}~\!$$ (existence of such function follows from \cite{Li84}). 
In turns out that $u^{*}\leq v$ on $\R^N$, 
see for instance Remark II.13 in  \cite{Li84}.  If 
 $u=u^{*}$ then we are done.  So assuming  
 $u\neq u^*$ we would like to derive a contradiction. By the theory of symmetrization (see Lieb and Loss \cite{LL}, Theorem 3.4),
we first obtain
$$
\int_{\Rn}|\De v|^p~\!dx=\int_{\Rn}|(-\De u)^{*}|^p~\!dx=\int_{\Rn}|\De u|^p~\!dx~.
$$
In addition, since we are assuming that $u^*\neq u$, then
$$
\int_{\Rn}|x|^{-2p}|u|^p~\!dx < \int_{\Rn}|x|^{-2p}(|u|^p)^*dx= \int_{\Rn}|x|^{-2p}|u^*|^p dx\leq \int_{\Rn}|x|^{-2p}|v|^pdx~\!.
$$
Thus we infer that
$$
\lambda\int_{\Rn}|x|^{-2p}|u|^p~\!dx\le\lambda\int_{\Rn}|x|^{-2p}|v|^p~\!dx~\!,
$$
and that the strict inequality holds if $\lambda>0$.
Similarly, we find
$$
\int_{\Rn}|x|^{-\beta}|u|^q~\!dx\le \int_{\Rn}|x|^{-\beta}|v|^q~\!dx~\!,
$$
and the strict inequality holds if $\beta>0$, that is, if $q<\starstar$.
In conclusion, since we are assuming that $\lambda$ and $\beta$ are not
contemporarily zero, we have that
$$
S_{q}(\la)\le
\frac{\displaystyle \int_{\Rn} \displaystyle\left[|\De v|^p-\la|x|^{-2p}|v|^p \right]dx}{
\left(\displaystyle \int_{\Rn}|x|^{-\beta}|v|^qdx\right)^{2/q}}\\
<
\frac{\displaystyle \int_{\Rn}[|\De u|^p-\la|x|^{-2p}|u|^p]dx}{
\left(\displaystyle \int_{\Rn}|x|^{-\beta}|u|^qdx\right)^{2/q}}=S_{q}(\la)~\!,
$$
a contradiction.
Therefore $u=u^{*}$, that is, $u$ is a nonnegative and radially symmetric decreasing function. 

\QED

As soon as $\lambda\to -\infty$, a braking symmetry phenomenon appears. In the next theorem we study 
the case $q<\starstar$, due to the
nonexistence result pointed out in the critical case $q=\starstar$,
$\lambda<0$. We cite \cite{DELT} , \cite{FelSch} for remarkable
breaking symmetry results for similar second-order equations in the case $p=2$. Also see \cite{BM} and \cite{C} for the similar type of results in the case of biharmonic equations.

\begin{Theorem}
If $\la<<0$ and $2\leq p<q<\starstar$ then $S_{q}(\la)<S_{q}^{\rm rad}(\la)$ and hence no extremal for $S_{q}(\la)$ is radially symmetric. 
\end{Theorem}

\proof
We already know that $S_{q}(\la)\leq S_{q}^{\rm rad}(\la)$. We will give 
an explicit condition
on $\lambda$ to have $S_{q}(\la)<S_{q}^{\rm rad}(\la)$.
Define 
$$
n(u)=\int_{\Rn}[|\De u|^p-\la|x|^{-2p}|u|^p]~\!dx~,\quad d(u)= \left(\int_{\Rn}|x|^{-\ba}|u|^q~\!dx\right)^{p/q}
$$ 
and $Q(u)=n(u)/d(u)$. 
Let
$u$ be a radially symmetric minimizer of $Q$ on $\mathcal D^{2,p}(\Rn)$. Our goal is to show that
$-\lambda$ can not be too large.  By homogeneity we can assume that $d(u)=1$. Thus 
$Q{'}(u)\cdot v=0$ and $ Q{''}(u)[v, v]\geq 0$ for all $v\in D^{2,p}(\Rn)$, that is,
 and
\begin{eqnarray}\label {c}
n{'}(u)\cdot v &=& Q(u) d{'}(u)\cdot v \nonumber\\
n{''}(u)[v, v] &\geq& Q(u) d{''}(u)[v, v]
\end{eqnarray}
for all $v\in D^{2,p}(\Rn)$.

 Let $\varphi_1\in H^{1}(\mathbb S^{N-1})$ be an Eigenfunction of Laplace-Beltrami operator on $\mathbb S^{N-1}$ corresponding to the smallest positive Eigenvalue. Thus
$$-\De_{\sigma}\varphi_1=(N-1)\varphi_1~, \   \  \frac{1}{|\mathbb S^{N-1}|}\int_{\mathbb S^{N-1}}|\varphi_1|^{2}~\!d\sigma=1~, \   \  \int_{\mathbb S^{N-1}}\varphi_1~\!d\sigma=0.$$
Now we set the test function $v$ as $v(x):=u(|x|)\varphi_1(\frac{x}{|x|})$. Therefore it turns out that
\begin{eqnarray}\label{d-v}
d{''}(u)[v, v] &=& p(p-q)\displaystyle\left(\int_{\Rn}|x|^{-\ba}|u|^{q-2}uv ~\!\right)^2
+p(q-1)\int_{\Rn}|x|^{-\ba}|u|^{q-2}v^2 \nonumber\\
&=& p(q-1)\int_{\Rn}|x|^{-\ba}|u|^{q}=p(q-1).
\end{eqnarray}
Also we see that $$|\De v|^2=|\De(u\varphi_1)|^2= |\De u-(N-1)|x|^{-2}u|^2\varphi_1^2$$
and therefore as $p\geq 2$, we have
\begin{eqnarray*}
n''(u)[v,v] &=& p(p-1)\displaystyle\int_{\Rn}\displaystyle\left(|\De u|^{p-2}|\De v|^2-\la|x|^{-2p}|u|^{p-2}|v|^2\right) dx\\
&=&p(p-1)\displaystyle\left[\int_{\Rn}|\De u|^{p-2}\big|\De u-(N-1)\frac{u}{|x|^{2}}\big|^2\varphi_1^2-\la|x|^{-2p}|u|^p\varphi_1^2 \right]  dx\\
&=&p(p-1)\int_{\Rn}\displaystyle\left[|\De u|^{p-2}\big|\De u-(N-1)\frac{u}{|x|^2}\big|^2-\la|x|^{-2p}|u|^p\right]dx\\
&=&p(p-1)\int_{\Rn}\displaystyle\left[|\De u|^p+\frac{(N-1)^2u^2|\De u|^{p-2}}{|x|^4}-\frac{2(N-1)|\De u|^{p-2}u\De u}{|x|^2}-\la\frac{|u|^p}{|x|^{2p}} \right]dx
\end{eqnarray*}
Since $p\geq 2$, using H\"older inequality and the fact that $d(u)=1$  we obtain,
\begin{eqnarray*}
n''(u)[v,v] &\leq& p(p-1)\displaystyle\left[n(u)+2(N-1)\displaystyle\left(\int_{\Rn}|\De u|^p dx\right)^\frac{p-1}{p}\displaystyle\left(\int_{\Rn}\frac{|u|^p}{|x|^{2p}}dx\right)^\frac{1}{p}\right]\\
&+&p(p-1)(N-1)^2\displaystyle\left(\int_{\Rn}|\De u|^p dx\right)^\frac{p-2}{p}\displaystyle\left(\int_{\Rn}\frac{|u|^p}{|x|^{2p}}dx\right)^\frac{2}{p}.
\end{eqnarray*}

Therefore from \eqref{c}, \eqref{d-v} and the definition of $Q(u)=n(u)$ we obtain
\begin{eqnarray*}
\nonumber
(q-p)Q(u) &\le& 2(p-1)(N-1)\displaystyle\left(\int_{\Rn}|\De u|^p dx\right)^\frac{p-1}{p}\displaystyle\left(\int_{\Rn}\frac{|u|^p}{|x|^{2p}}dx\right)^\frac{1}{p}\\
&+&(p-1)(N-1)^2 \displaystyle\left (\int_{\Rn}|\De u|^p dx\right)^\frac{p-2}{p}\displaystyle\left(\int_{\Rn}\frac{|u|^p}{|x|^{2p}} dx\right)^\frac{2}{p}.
\end{eqnarray*}
Thus we have,
\begin{eqnarray*}
\nonumber
(q-p)\int_{\Rn}|\De|^p dx &\le& \la(q-p)\int_{\Rn}\frac{|u|^p}{|x|^{2p}}+ 2(p-1)(N-1)\displaystyle\left(\int_{\Rn}|\De u|^p dx\right)^\frac{p-1}{p}\displaystyle\left(\int_{\Rn}\frac{|u|^p}{|x|^{2p}}dx\right)^\frac{1}{p}\\
&+&(p-1)(N-1)^2 \displaystyle\left(\int_{\Rn}|\De u|^p dx\right)^\frac{p-2}{p}\displaystyle\left(\int_{\Rn}\frac{|u|^p}{|x|^{2p}} dx\right)^\frac{2}{p}.
\end{eqnarray*}

In particular, the quantity
$$
X:=\left(\frac{\displaystyle\int_{\Rn}|\De u|^p ~\!}
{\displaystyle\int_{\Rn}|x|^{-2p}|u|^p~\!}\right)^{1/p}
$$
satisfies the inequality
$$
(q-p)X^p\leq\la(q-p)+2(N-1)(p-1)X^{p-1}+(p-1)(N-1)^2X^{p-2} ~\!,
$$ which implies
$$
\lambda\ge \min_{t\in\R} f(t),$$ where
$$f(t)=\left(~\!t^p-\frac{2(N-1)(p-1)}{q-p}t^{p-1}-\frac{(p-1)(N-1)^2}{q-p}t^{p-2}\right).
$$
Define
 $$t_0=\frac{\gamma_1(p-1)+\sqrt{\gamma^2_1(p-1)^2+4\gamma_2p(p-2)}}{2p}$$ where 
 $\gamma_1=\frac{2(p-1)(N-1)}{q-p}$ and $\gamma_2=\frac{(p-1)(N-1)^2}{q-p}$. By a straight forward calculation we obtain $\gamma\geq f(t_0)$.

Therefore no extremal for $S_{q}(\la)$ is radially symmetric and break of symmetry occurs if $\la<f(t_0)$.
\QED

\begin{Remark}
If $q$ is close enough to $\starstar$ then one can obtain a better estimate
on the breaking symmetry parameter $\lambda$ by arguing as follows.
Notice that 
$X> \gamma_{N,p}$ by the Rellich inequality  (\ref{Rellich}). Thus, if
$\gamma_{N,p}\ge t_0$ that is, if
$$
p+\frac{p(N-1)\displaystyle\left[p^2(N-1)(p-2)+2N(p-1)^2(N-2p)\right]}{N^2(p-1)(N-2p)^2}\le q\le\starstar~\!,
$$
then the radial solution $u$ does not achieve $S_q(\lambda)$ 
unless
$$
\lambda> \min_{t\ge \gamma_{N,p}} f(t)
=\gamma_{N,p}^p-\gamma_1(\gamma_{N,p})^{p-1}-\gamma_2(\gamma_{N,p})^{p-2}
$$ Conversely, if 
$$
\la\le \gamma_{N,p}^p-\gamma_1(\gamma_{N,p})^{p-1}-\gamma_2(\gamma_{N,p})^{p-2}
$$
then break of symmetry occurs.
\end{Remark}

\small
\noindent
{\bf Acknowledgements.} 
The  Author wishes to thank Prof. Roberta Musina for having suggested the references \cite{Miti} , \cite{Mu1}, \cite{Mu2}  and for her useful comments. the author also wishes to thank the anonymous referee for careful reading the manuscript and for his/her many valuable suggestions. This research was supported by the fund of Australian Research Council (ARC).

\label{References}


\begin{thebibliography}{XX}


\bibitem{AlDo} {Alves, C. O. ;  do \'{O}, Jo\~{a}o Marcos}, {Positive solutions of a fourth-order semilinear problem involving critical growth}, \textit{Adv. Nonlinear Stud.}, {\textbf 2}, (2002), no 4, 437--458.


\bibitem{BM} {Bhakta, Mousomi; Musina, Roberta}, {Entire solutions for a class of variational problems involving the biharmonic operator and Rellich potentials}, Nonlinear Anal. {\textbf 75} (2012), no. 9, 3836--3848.

\bibitem{CKN}
{Caffarelli, L.; Kohn, R.; Nirenberg, L.},
{First Order Interpolation Inequalities with Weights},
\textit{Compositio Math.} \textbf{53} (1984), 259--275.

\bibitem{C}{Caldiroli, P.}, {Radial and non radial ground states for a class of dilation invariant fourth order semilinear elliptic equations on  $\Rn$}, \textit{Communications on Pure and Applied Analysis}	, Volume {\textbf 13}, (2014), no 2, 811--821.

\bibitem{CM2}
{Caldiroli, P.; Musina, R.}, 
Caffarelli-Kohn-Nirenberg type inequalities for the weighted biharmonic operator in cones,
\textit{Milan J. Math.},  {\bf 79} (2011), no 2, 657--687.


\bibitem{CatWan}
{Catrina, F.; Wang, Z.-Q.}, 
{On the Caffarelli-Kohn-Nirenberg inequalities: sharp constants, existence (and nonexistence), and symmetry of extremal functions},
\textit{Comm. Pure Appl. Math.} \textbf{54} (2001), 229--258.


\bibitem{Cof}
{Coffman, Ch. D.},
{On the structure of solutions to $\De^2u = \la u$ which satisfy the clamped plate conditions on a right angle}, \textit{SIAM J. Math. Anal} \textbf{13} (1982), 746--757. 

\bibitem{DH}{Davies, E. B; Hinz, A. M.}  {Explicit constants for Rellich inequalities in $L^p(\Omega)$}, \textit{Math. Z.} \textbf {227}, (1998) 511--523.

\bibitem{DELT} {Dolbeault, J.; Esteban, M.; Loss, M.; Tarantello, G.} {On the symmetry of extremals for the Caffarelli-Kohn-Nirenberg inequalities}, \textit{Adv. Nonlinear Stud.} \textbf{9}, (2009), 713--726.

\bibitem{FelSch}
{Felli, V.; Schneider, M.}, 
{Perturbation results of critical elliptic equations of Caffarelli-Kohn-Nirenberg type}, 
\textit{J. Diff. Eq.}
\textbf{191} (2003), 121--142. 

\bibitem{GG}{Gazzola, F.; Grunau, H. C.}, {Radial entire solutions for supercritical biharmonic equations,}
\textit{Math. Ann.} \textbf{334} (2006), 905-936.

\bibitem{HV} {Hulshof, J.;  Van der Vorst, R.C.A.M.}, {Asymptotic behaviour of ground states}, \textit{Proc. Amer.
Math. Soc.} \textbf{124}, (1996), no. 8, 2423--2431.


\bibitem{LL} {Lieb, H. E., Loss, M.}, {Analysis, Second Edition}, \textit{Graduate Studies in Mathematics} Vol 14, American Mathematical Society, 2001.

\bibitem{Li84} Lions, P. L., The concentration-compactness principle in the calculus of variations.The locally compact case. II. \textit{Ann. Inst. H. Poincar\'e Anal. Non
Lin\'eaire} \textbf{1} (1984),  223--283.

\bibitem{Miti}{Mitidieri, E.}, {A simple approach to Hardy inequalities}, \textit{Mat. Zametki} \textbf{ 67}, (2000), no. 4, 563--572; translation in Math. Notes { 67} (2000), no.~3-4, 479--486.

\bibitem{Mu1}{Musina, Roberta}, {Optimal Rellich-Sobolev constants and their extremals},  \textit{Differential Integral Equations} {\bf 27} (2014), no. 5-6, 579--600.

\bibitem{Mu2}{Musina, Roberta}, {Weighted Sobolev spaces of radially symmetric functions}, \textit{Ann. Mat. Pura Appl.}, (2013)  arXiv:1206.6957.

\bibitem{Rel54}
{Rellich, F.}, 
Halbbeschr\"ankte Differentialoperatoren h\"oherer Ordnung. 
In: J.C.H. Gerretsen, J. de Groot (Eds.): Proceedings of the International Congress of Mathematicians 1956, Volume III 
(pp. 243--250) Groningen: Noordhoff 1956.

\bibitem{Rel69}
{Rellich, F.},
{Perturbation theory of eigenvalue problems}, \textit{New York: Courant Institute of Mathematical Sciences, New York University}, 
 (1954).
 
\end{thebibliography}
\end{document}